\documentclass[12pt]{amsart}

\usepackage{amsmath,amsthm,amsfonts,latexsym,epsfig}
\newtheorem{theorem}{Theorem}[section]
\newtheorem{lemma}[theorem]{Lemma}
\newtheorem{proposition}[theorem]{Proposition}
\newtheorem{corollary}[theorem]{Corollary}
\newtheorem{conjecture}[theorem]{Conjecture}

\title{\hbox{Degree Bounds in Quantum Schubert Calculus}}
\author{Alexander Yong
}
\date{\today} 
\address{Department of Mathematics, University of Michigan, Ann Arbor, Michigan
48109
}

\email{ayong@umich.edu}

\begin{document}
\input amssym.def
\input amssym.tex
\begin{abstract}
Fulton and Woodward \cite{Fulton_Woodward} have recently identified the smallest
degree of $q$ that appears in the expansion of the product of two Schubert classes 
in the (small) quantum cohomology ring of a Grassmannian. 
We present a combinatorial proof of this result, and provide an
alternative characterization of this smallest degree
in terms of the rim hook formula~\cite{BCF} for the quantum
product.  
\end{abstract}

\maketitle

\section{Introduction and Main Results}

        Let $X=Gr(l,{\mathbb C}^{n})$ be the Grassmannian of $l$-dimensional subspaces in
${\mathbb C}^{n}$. The classical cohomology
ring ${\rm H}^{*}(X,{\mathbb Z})$ has an additive basis of {\em Schubert classes} 
$\{\sigma_{\lambda}\}$, indexed by the Young diagrams $\lambda$ (identified with the corresponding partitions)
contained in the $l\times k$ rectangle, where $k=n-l$ (we denote this by $\lambda\subseteq l\times k$). 
The product of two Schubert classes in ${\rm H}^{*}(X,{\mathbb Z})$ is given by
\begin{equation}
\label{LR_exp}
\sigma_{\lambda}\cdot\sigma_{\mu}=\sum_{\nu\subseteq l\times k}c_{\lambda,\mu}^{\nu}\sigma_{\nu},
\end{equation} 
where $c_{\lambda,\mu}^{\nu}$ is the Littlewood-Richardson coefficient (see, e.g., \cite{YT,EC2}).

        The (small) quantum cohomology ring ${\rm QH}^{*}(X)$ is a certain deformation
of ${\rm H}^{*}(X,{\mathbb Z})$ that has been extensively studied in recent years; see, 
e.g., 
\cite{Bertram,Fulton_P,Sottile} and references therein. 
This ring is canonically isomorphic to the
Verlinde algebra of $\mathfrak{sl}_{l}$ at level $k$ (see, e.g.,
\cite{BCF,Witten}); consequently, the results presented
here can be reformulated in representation theoretic language.
 
        The additive structure of ${\rm QH}^{*}(X)$ is essentially
the same as that of ${\rm H}^{*}(X,{\mathbb Z})$: the Schubert classes form a basis 
of ${\rm QH}^{*}(X)$ as a free module over ${\mathbb Z}[q]$, where $q$ is an indeterminate.
The multiplicative structure of ${\rm QH}^{*}(X)$ is defined by

\begin{equation}
\label{ring_collected}
\sigma_{\lambda}*\sigma_{\mu}=\sum_{\nu\subseteq l\times k}\sum_{d\geq 0}q^{d}
\langle\lambda,\mu,\nu^{\vee}\rangle_{d}\sigma_{\nu},
\end{equation}
where $\nu^{\vee}=(k-\nu_{l},\ldots,k-\nu_{1})$ is the complement of $\nu$ in the $l\times k$
rectangle, and  $\langle \lambda,\mu,\nu^{\vee}\rangle_{d}$
is a three-point, genus-zero Gromov-Witten invariant of $X$. 
Note that we use ``$*$'' to distinguish the quantum multiplication from the product 
in the classical
cohomology ring. Setting $q=0$ recovers (\ref{LR_exp})
because $c_{\lambda,\mu}^{\nu}=\langle\lambda,\mu,\nu^{\vee}\rangle_{0}$. 

        Bertram, Ciocan-Fontanine and Fulton \cite{BCF} have given a combinatorial rule to 
compute $\sigma_{\lambda}*\sigma_{\mu}$ and thus the Gromov-Witten invariants
$\langle\lambda,\mu,\nu^{\vee}\rangle_{d}$.  
To describe this rule, we need some terminology and
notation. An $n${\em-rim hook} of a Young diagram $\lambda$ is a connected subset of $n$ boxes 
of $\lambda$ that does not contain a $2\times 2$ square. The {\em width} of an $n$-rim hook is 
the number of columns it occupies. An $n$-rim hook is {\em legal} if
removing it from a Young diagram gives a valid Young diagram (see Figure 1). 
For a partition $\rho$, we define its $n$-core, denoted 
${\rm core}_{n}(\rho)$, to be the partition corresponding to the Young diagram  
obtained by repeatedly removing legal $n$-rim hooks from $\rho$ until  
further removals are not possible. It is well known (see, e.g., \cite{James_Kerber}) 
that this procedure defines ${\rm core}_{n}(\rho)$ uniquely. 
Let $r_{n}(\rho)=\frac{\mid\rho\mid -\mid {\rm core}_{n}(\rho)\mid}{n}$ 
be the number of $n$-rim hooks removed in this process, 
and set $\epsilon(\rho)=(-1)^{\sum(k-width(R_{i}))}$, where 
$R_{1},\ldots,R_{r_{n}(\rho)}$ are these $n$-rim hooks. With this notation, 
the rule obtained in \cite{BCF} is as follows:
\begin{equation}
\label{ring_structure}
\sigma_{\lambda}*\sigma_{\mu}=\sum_{\nu\subseteq l\times k}
\sum_{\begin{array}{c}
\scriptstyle{c_{\lambda,\mu}^{\rho}\neq 0} \\
\scriptstyle{\rho_{1}\leq k} \\
\scriptstyle{{\rm core}_{n}(\rho)=\nu}
\end{array}
} 
q^{r_{n}(\rho)}\epsilon(\rho)c_{\lambda,\mu}^{\rho}\sigma_{\nu}.
\end{equation}

The Gromov-Witten invariants $\langle \lambda,\mu,\nu^{\vee}\rangle$ are known
to be nonnegative; thus the coefficient of each term $q^{d}\sigma_{\nu}$ in the right-hand 
side of (\ref{ring_structure}) is nonnegative.
It remains an open problem to find a direct combinatorial 
proof of this statement; such a proof would provide a 
generalization of the Littlewood-Richardson rule to the
quantum Schubert calculus of the Grassmannian.

\begin{figure}[t]
\centering
\epsfig{file=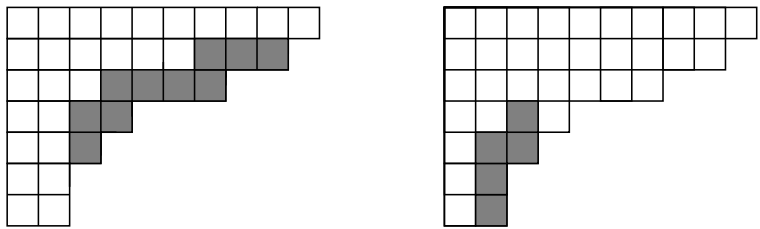,height=3cm}
\caption{Legal and illegal rim hooks for $\lambda=(10,9,7,4,3,2,2)$}
\end{figure}

        Let ${\rm d_{min}}$ denote the smallest
degree of $q$ such that $q^{\rm{d_{min}}}$ appears in (\ref{ring_collected}) 
with nonzero coefficient. The
following theorem provides an affirmative answer to a question posed by
A.~Ram (private communication).\pagebreak

\begin{theorem}
\label{main_thm}
Let $\lambda,\mu\subseteq l\times k$.
Among all $\rho$ with $\rho_{1}\leq k$ and \linebreak
$c_{\lambda,\mu}^{\rho}\neq 0$, pick those with the largest
n-core (equivalently, the smallest value of~$r_{n}(\rho)$). 
Then their contributions to {\rm (\ref{ring_structure})} do not cancel each other out.
In other words, ${\rm d_{min}}$ equals the smallest degree of $q$ that 
appears in the right-hand side of
{\rm (\ref{ring_structure})} (before cancellations).
\end{theorem}

        Our proof of Theorem~\ref{main_thm} will also prove a recent result of Fulton and 
Woodward (see Theorem~1.2 below).
For partitions $\lambda,\mu\subseteq l\times k$, place $\lambda$ against the upper 
left corner of the rectangle.
Then rotate $\mu$ 180 degrees and place it in the lower right corner (see Figure~2). We will refer to 
${\rm rotate}(\mu)$ as the resulting subshape of $l\times k$. 
Let ${\mathfrak d}$ be the side length of the largest square that fits inside $\lambda\cap{\rm rotate}(\mu)$.
The following theorem was conjectured by Fulton, and later proved by Fulton and Woodward 
\cite{Fulton_Woodward}, 
using moduli spaces.\footnote{Actually, Fulton and Woodward proved
a more general result that applies to any homogeneous space $X=G/P$, where $G$ is a simply
connected complex semisimple Lie group and $P$ its parabolic subgroup.} 


\begin{theorem}
\label{second_main_thm}
{\rm \cite{Fulton_Woodward}}
${\rm d_{min}}={\mathfrak d}$.
\end{theorem}

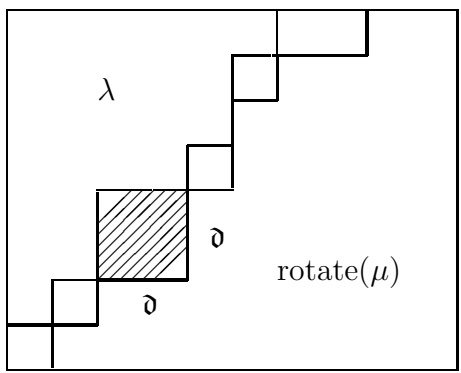
\begin{figure}[h]
\setlength{\unitlength}{.85pt} 
\centering
\begin{picture}(220,180)
\put(15,15){\framebox(200,160)}
\thicklines
\put(15,35){\line(1,0){40}}
\put(55,35){\line(0,1){20}}
\put(55,55){\line(1,0){40}}
\put(95,55){\line(0,1){60}}
\put(95,115){\line(1,0){20}}
\put(115,115){\line(0,1){20}}
\put(115,135){\line(1,0){20}}
\put(135,135){\line(0,1){20}}
\put(135,155){\line(1,0){40}}
\put(175,155){\line(0,1){20}}
\put(55,135){$\lambda$}
\thinlines
\put(55,55){\line(1,1){40}}
\put(55,60){\line(1,1){35}}
\put(55,65){\line(1,1){30}}
\put(55,70){\line(1,1){25}}
\put(55,75){\line(1,1){20}}
\put(55,80){\line(1,1){15}}
\put(55,85){\line(1,1){10}}
\put(55,90){\line(1,1){5}}

\put(60,55){\line(1,1){35}}
\put(65,55){\line(1,1){30}}
\put(70,55){\line(1,1){25}}
\put(75,55){\line(1,1){20}}
\put(80,55){\line(1,1){15}}
\put(85,55){\line(1,1){10}}
\put(95,55){\line(1,1){5}}
\put(35,15){\dashbox{.5}(0,40)}
\put(35,55){\dashbox{.5}(20,0)}
\put(55,55){\dashbox{.5}(0,40)}
\put(55,95){\dashbox{.5}(60,0)}
\put(115,95){\dashbox{.5}(0,60)}
\put(115,155){\dashbox{.5}(20,0)}
\put(135,155){\dashbox{.5}(0,20)}
\put(135,55){${\rm rotate}(\mu)$}
\put(75,40){${\mathfrak d}$}
\put(105,70){${\mathfrak d}$}
\end{picture}
\caption{$\lambda$, ${\rm rotate}(\mu)$ and ${\mathfrak d}\times {\mathfrak d}$}
\end{figure}

        An alternative proof of Theorem~\ref{second_main_thm} was later given by
A.~Buch \cite{Buch_preprint}, using an elegant geometric argument combined with
combinatorics. We will utilize the combinatorial part of Buch's proof below.

        Let us now describe our proof of Theorems~\ref{main_thm} and \ref{second_main_thm}. 
The proof is entirely combinatorial once 
the nonnegativity of the Gromov-Witten invariants and the fact that ${\rm QH}^{*}(X)$ is an
associative ring \cite{Kontsevich-Manin,Ruan-Tian} are granted. Our argument is based on
the following result.   


\begin{theorem}
\label{first_main_thm}
Let $a\times A$, $b\times B$ and $c\times C$ be rectangular Young
diagrams contained in $l\times k$.
Then the following are equivalent:
\begin{itemize}
\item[{\rm (I)}] there exist Young diagrams $\lambda,\mu,\nu\subseteq l\times k$ 
containing $a\times A$, $b\times B$ and $c\times C$, respectively, such that 
$\sigma_{\lambda}\cdot\sigma_{\mu}\cdot\sigma_{\nu}\neq 0$;
\item[{\rm (II)}] $\sigma_{a\times A}\cdot\sigma_{b\times
    B}\cdot\sigma_{c\times C}\neq 0$;
\item[{\rm (III)}] all of the five conditions below hold:
\begin{itemize}
\item[{\rm (i)}] $a+b\leq l$ or $A+B\leq k$;
\item[{\rm (ii)}] $a+c\leq l$ or $A+C\leq k$;
\item[{\rm (iii)}] $b+c\leq l$ or $B+C\leq k$;
\item[{\rm (iv)}] $a+b+c\leq l$ or $A+B+C\leq 2k$;
\item[{\rm (v)}] $a+b+c\leq 2l$ or $A+B+C\leq k$.
\end{itemize}
\end{itemize}
\end{theorem}
                                             

        A combinatorial proof of
        Theorem~\ref{first_main_thm} is given in Section~2.
S.~Fomin (private communication) has observed that 
the necessity of the conditions (i)-(v) can be derived 
from the Horn inequalities describing the ``Klyachko cone'' 
$\{(\lambda,\mu,\nu)\mid c_{\lambda,\mu}^{\nu}\neq 0\}$; specifically see
\cite[(11)]{Fulton_BAMS}.

\begin{corollary}
\label{Cor_sect2}
Let $\lambda,\mu,\rho$ be partitions such that $\lambda,\mu\subseteq l\times k$
and $\rho_{1}\leq k$. If $c_{\lambda,\mu}^{\rho}\neq 0$, 
then $(l+{\mathfrak d})\times {\mathfrak d}\subseteq \rho$. 
\end{corollary}  

\begin{proof}
Suppose $c_{\lambda,\mu}^{\rho}\neq 0$ but 
$(l+{\mathfrak d})\times {\mathfrak d}\not\subseteq \rho$. Choose any positive
integer $L\geq l$ such that $\rho\subseteq L\times k$. Then we have $\sigma_{\lambda}\cdot
\sigma_{\mu}\cdot\sigma_{\rho^{\vee}}\neq 0$. Notice that $(L-(l+{\mathfrak d})+1)\times (k-{\mathfrak d}+1)\subseteq \rho^{\vee}$. By 
the definition of ${\mathfrak d}$, there exist rectangles $a\times A\subseteq\lambda$ and
$b\times B\subseteq \mu$ such that $a+b\geq l+{\mathfrak d}$ and
$A+B\geq k+{\mathfrak d}$. 
Set $\nu=\rho^{\vee}, c=L-(l+{\mathfrak d})+1$
and $C=k-{\mathfrak d}+1$, then it is
easy to check that both corresponding inequalities (iv) 
are violated, a contradiction of Theorem~\ref{first_main_thm} (see Figure~3 below). 
\end{proof}

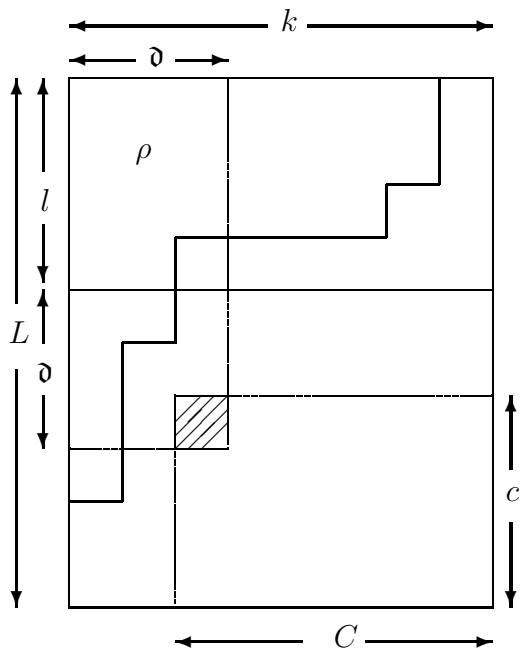
\begin{figure}[h]
\label{Corollary_figure}
\centering
\begin{picture}(200,240)(-15,0)
\put(25,15){\framebox(160,200)}
\thinlines
\put(65,15){\dashbox{.5}(0,80)}
\put(65,95){\dashbox{.5}(120,0)}
\put(25,75){\dashbox{.5}(60,0)}
\put(85,75){\dashbox{.5}(0,140)}
\put(25,135){\line(1,0){160}}
\thicklines
\put(25,55){\line(1,0){20}}
\put(45,55){\line(0,1){60}}
\put(45,115){\line(1,0){20}}
\put(65,115){\line(0,1){40}}
\put(65,155){\line(1,0){80}}
\put(145,155){\line(0,1){20}}
\put(145,175){\line(1,0){20}}
\put(165,175){\line(0,1){40}}
\put(13,100){${\mathfrak d}$}
\put(15,95){\vector(0,-1){20}}
\put(15,115){\vector(0,1){20}}
\put(14,165){$l$}
\put(15,160){\vector(0,-1){22}}
\put(15,180){\vector(0,1){35}}
\put(105,233){$k$}
\put(115,235){\vector(1,0){70}}
\put(100,235){\vector(-1,0){75}}
\put(55,220){${\mathfrak d}$}
\put(50,222){\vector(-1,0){25}}
\put(65,222){\vector(1,0){20}}
\put(125,0){$C$}
\put(137,2){\vector(1,0){48}}
\put(113,2){\vector(-1,0){48}}
\put(190,55){$c$}
\put(192,67){\vector(0,1){28}}
\put(192,47){\vector(0,-1){32}}
\put(50,185){$\rho$}

\put(5,130){\vector(0,1){85}}
\put(5,105){\vector(0,-1){90}}
\put(2,115){$L$}

\thinlines
\put(65,75){\line(1,1){20}}
\put(65,80){\line(1,1){15}}
\put(65,85){\line(1,1){10}}
\put(65,90){\line(1,1){5}}
\put(70,75){\line(1,1){15}}
\put(75,75){\line(1,1){10}}
\put(80,75){\line(1,1){5}}
\end{picture}
\caption{Proof of Corollary \ref{Cor_sect2}.}
\end{figure}


\medskip
\noindent \emph{Proof of Theorems~\ref{main_thm} and
  \ref{second_main_thm}.} 
Let $\rho$ be a partition whose first part is
at most $k$ and whose $n$-core is contained in $l\times k$. 
If $c_{\lambda,\mu}^{\rho}\neq 0$, then Corollary~\ref{Cor_sect2} gives
$(l+{\mathfrak d})\times {\mathfrak d}\subseteq \rho$, implying that 
at least ${\mathfrak d}$ rim hook removals are needed to
obtain ${\rm core}_{n}(\rho)$. Hence the smallest degree of $q$ that occurs
in (\ref{ring_structure}), and therefore ${\rm d_{min}}$,
is at least ${\mathfrak d}$.

        To complete the proofs, it 
remains to show that ${\rm d_{min}}$ and thus the smallest degree of $q$ that appears in 
(\ref{ring_structure}) is at most ${\mathfrak d}$. To this end,
we borrow an argument of A. Buch \cite{Buch_preprint}, used in his own proof
of Theorem~\ref{second_main_thm}. We reproduce his argument below:
 
        Let $\widetilde{\lambda}$ and $\overline{\lambda}$ be the partitions obtained by
removing the leftmost ${\mathfrak d}$ columns and the top ${\mathfrak d}$ rows of $\lambda$,
respectively. Also, let $\widetilde{\overline \lambda}$ be the partition obtained by
removing both the leftmost ${\mathfrak d}$ columns and the top ${\mathfrak d}$ rows of $\lambda$. 
Set $\alpha=(k+{\mathfrak d}-\lambda_{\mathfrak d},\ldots,k+{\mathfrak d}-\lambda_{1})$ and let \linebreak
$\beta_{i}=\max({\mathfrak d}-\lambda_{l+1-i},0)$ for $1\leq i\leq l$. In other words, 
$\alpha$ is the complement of the bottom ${\mathfrak d}$ rows of ${\rm rotate}(\lambda)$ in 
$l\times (k+{\mathfrak d})$, and $\beta$ is the complement of ${\rm rotate}(\lambda)$ in the
rightmost ${\mathfrak d}$ columns (see Figure 4 below).

        It follows from the Littlewood-Richardson rule that the expansion of  
$\sigma_{\lambda}*\sigma_{\beta}$ contains the class
$\sigma_{({\mathfrak d}^{l})+\widetilde{\lambda}}
=\sigma_{({\mathfrak d}^{l})}*\sigma_{\widetilde{\lambda}}$.
It also follows that $\sigma_{\widetilde{\lambda}}*\sigma_{\alpha}$ contains
$\sigma_{(k^{{\mathfrak d}}),\widetilde{\overline{\lambda}}}=
\sigma_{(k^{{\mathfrak d}})}*\sigma_{\widetilde{\overline{\lambda}}}$.  
It is not hard to check directly from (\ref{ring_structure}) that 
$\sigma_{({\mathfrak d}^{l})}*\sigma_{(k^{{\mathfrak d}})}=q^{\mathfrak d}$.
By the nonnegativity of the Gromov-Witten invariants and the associativity
of ${\rm QH}^{*}(X)$,
$\sigma_{\lambda}*\sigma_{\beta}*\sigma_{\alpha}$ contains the product
$\sigma_{({\mathfrak d}^{l})}*\sigma_{(k^{{\mathfrak d}})}*
\sigma_{\widetilde{\overline \lambda}}=q^{{\mathfrak d}}
\sigma_{\widetilde{\overline \lambda}}$. 
Note that $\widetilde{\overline \lambda}\cap {\rm rotate}(\mu)=\emptyset$, which is well known
to be equivalent to 
$\sigma_{\widetilde{\overline \lambda}}\cdot\sigma_{\mu}\neq 0$. 
Thus we conclude that
$\sigma_{\lambda}*\sigma_{\mu}*\sigma_{\alpha}*\sigma_{\beta}$ contains 
$q^{\mathfrak d}$ times some Schubert class. Therefore, the product
$\sigma_{\lambda}*\sigma_{\mu}$ must have a term involving a degree of
$q$ less than or equal to ${\mathfrak d}$, and we are done.\qed 


\begin{figure}[h]
\label{Buch_fig}
\setlength{\unitlength}{.9pt} 
\centering
\begin{picture}(200,200)(-70,0)
\put(-160,130){$A\cup B\cup C\cup D={\rm rotate}(\lambda)$}
\put(-160,110){$A\cup B={\rm rotate}({\overline \lambda})$}
\put(-160,90){$A\cup C={\rm rotate}(\widetilde{\lambda})$}
\put(-160,70){$A={\rm rotate}(\widetilde{\overline \lambda})$}
\put(-160,50){$D={\rm rotate}({\mathfrak d}\times {\mathfrak d})$}
\put(15,15){\framebox(200,160)}
\thinlines
\put(55,15){\line(0,1){160}}
\put(175,14){\line(0,1){160}}
\put(15,55){\line(1,0){200}}
\thicklines
\put(75,15){\line(0,1){20}}
\put(75,35){\line(1,0){20}}
\put(95,35){\line(0,1){40}}
\put(95,75){\line(1,0){40}}
\put(135,75){\line(0,1){40}}
\put(135,115){\line(1,0){20}}
\put(155,115){\line(0,1){20}}
\put(155,135){\line(1,0){40}}
\put(195,135){\line(0,1){20}}
\put(195,155){\line(1,0){20}}
\put(5,90){$l$}
\put(110,180){$k$}
\put(32,2){${\mathfrak d}$}
\put(32,32){$\alpha$}
\put(193,2){${\mathfrak d}$}
\put(220,32){${\mathfrak d}$}
\put(193,162){$\beta$}
\put(145,75){$A$}
\put(145,32){$C$}
\put(193,75){$B$}
\put(193,32){$D$}
\end{picture}
\caption{Proof of Theorems~\ref{main_thm} and \ref{second_main_thm}}
\end{figure}
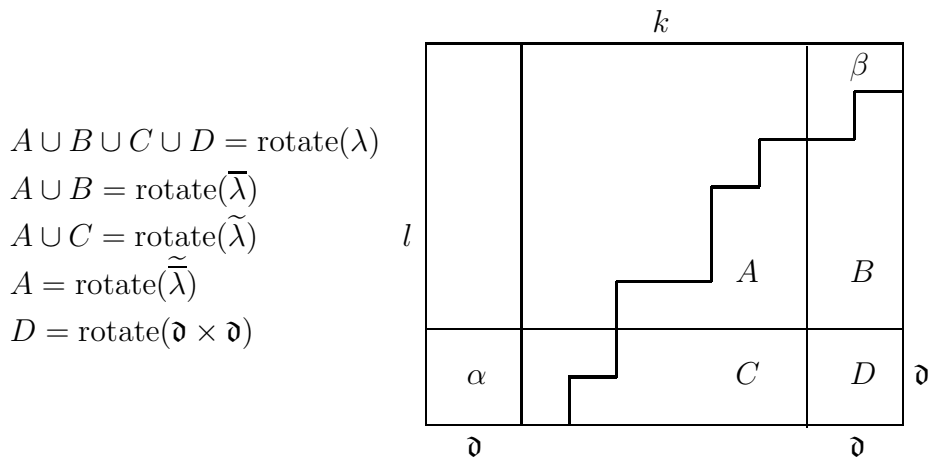

\section{Proof of Theorem~\ref{first_main_thm}}

        For $\alpha,\beta\subseteq l\times k$, define
$\varkappa(\alpha,\beta)$ to be the intersection of all Young
diagrams $\gamma\subseteq l\times k$ such that
$c_{\alpha,\beta}^{\gamma}\neq 0$. 
In the event the set of such~$\gamma$ is empty, we set
$\varkappa(\alpha,\beta)=l\times k$ by convention. 
We will make use of the following
observation, which follows immediately from this definition:

\begin{lemma}
\label{small_lemma}
Let $\alpha,\beta, m\times M\subseteq l\times k$. 
Then $\sigma_{\alpha}\cdot\sigma_{\beta}\cdot\sigma_{m\times M}\neq 0$
if and only if $\varkappa(\alpha,\beta)\cap{\rm rotate}(m\times M)=\emptyset$.
\end{lemma}

        Our proof of Theorem~\ref{first_main_thm} also uses the
        following two lemmas.\footnote{Although we will not need it here, 
this Lemma~2.2 generalizes to a statement about the cohomology of any homogeneous space $G/P$. 
The proof is analogous to the one given above.}

\begin{lemma}
\label{trick_lemma}
Let $\alpha,\beta,\lambda,\mu\subseteq l\times k$ be such that 
$\alpha\subseteq \lambda$ and
$\beta\subseteq \mu$. Then $\varkappa(\alpha,\beta)\subseteq \varkappa(\lambda,\mu)$.
\end{lemma}

\begin{proof}
It suffices to prove the case when $\lambda=\alpha$ and 
$\mu\setminus\beta$ is a single box. Since every term of 
$(\sigma_{\alpha}\cdot\sigma_{\beta})\cdot \sigma_{1}$
is indexed by a partition containing $\varkappa(\alpha,\beta)$, the same is true for
$\sigma_{\alpha}\cdot(\sigma_{\beta}\cdot\sigma_{1})=\sigma_{\lambda}\cdot
(\sigma_{\mu}+\mbox{nonnegative terms})$. The claim follows. 
\end{proof}

The next lemma is proved by a straightforward application of the
Littlewood-Richardson Rule. Details are left to the reader.

\pagebreak
\begin{lemma}
\label{kappa_calc_lemma}
Let $m\times M,n\times N\subseteq l\times k$. 
If $\sigma_{m\times M}\cdot\sigma_{n\times N}=0$, then
$\varkappa(m\times M,n\times N)=l\times k$.
Otherwise,
\begin{eqnarray} \nonumber
\varkappa(m\times M,n\times N) &=& (m\times M)\cup (n\times N)\\ \nonumber
& & \cup ((m+n)\times (M+N-k)) \\ \nonumber
 & & \cup ((m+n-l)\times (M+N)), \nonumber
\end{eqnarray}
where each of the last two rectangles in the right-hand side is
understood to be empty if one of its dimensions is negative,
or if it does not fit inside $l\times k$.
\end{lemma}

\noindent \emph{Proof of Theorem~\ref{first_main_thm}.} Suppose (I)
holds. Then there exists $\rho\subseteq
l\times k$ with $c_{\lambda,\mu}^{\rho}\neq 0$ and $\rho\cap
{\rm rotate}(\nu)=\emptyset$. Hence $\varkappa(\lambda,\mu)\cap
{\rm rotate}(\nu)=\emptyset$ which by Lemma~\ref{trick_lemma} implies
\begin{equation}
\label{star_eqn}
\varkappa(a\times A,b\times B)\cap {\rm rotate}(c\times C)=\emptyset.
\end{equation}
By Lemma~\ref{small_lemma}, (\ref{star_eqn}) is equivalent to (II),
which itself trivially implies~(I). Lastly, it follows from
Lemma~\ref{kappa_calc_lemma} that (\ref{star_eqn}) is equivalent to (III).\qed

\section{An Upper Bound, Open Problems and Conjectures}
        
        It is an open problem to determine ${\rm d_{max}}$,
the largest degree of $q$ such that $q^{{\rm d_{max}}}$ appears in
(\ref{ring_collected}) with nonzero coefficient.
In particular,  ${\rm d_{max}}$ is not equal to
the largest degree of $q$ appearing in (\ref{ring_structure}). 
Next, we present a simple upper bound for ${\rm d_{max}}$. For a Young diagram
$\alpha$, let ${\rm diag}(\alpha)$ be the size of the largest square contained in $\alpha$.

\begin{proposition}
\label{fourth_main_thm}
${\rm d_{max}}\leq \min({\rm diag}(\lambda),{\rm diag}(\mu)).$
\end{proposition}

        We note that for a large class of pairs $(\lambda,\mu)$, this inequality is sharper
than the obvious upper bound ${\rm d_{max}}\leq \frac{\mid\lambda\mid+\mid\mu\mid}{n}$.

\medskip
\noindent
\emph{Proof of Proposition~\ref{fourth_main_thm}} Let $\lambda$ be a
Young diagram and \linebreak $\lambda=(\alpha_{1},\ldots\alpha_{t}\mid \beta_{1},\ldots,\beta_{t})$
its {\em Frobenius notation}; the following is
a well known identity for Schur functions (see, e.g. \cite{Macdonald}):
\begin{equation}
\label{Giambelli}
s_{\lambda}=\det(s_{(\alpha_{i}\mid\beta_{j})})_{1\leq i,j\leq t}.
\end{equation}
From (\ref{Giambelli}) and our description of ${\rm QH}^{*}(X)$ it follows that 
\begin{equation}
\label{thm_2_two}
\sigma_{\lambda}=\det(\sigma_{(\alpha_{i}\mid\beta_{j})})_{1\leq i,j\leq t}.
\end{equation}
It is easy to check from the Littlewood-Richardson rule that
for any $(\alpha\mid\beta),\lambda\subseteq l\times k$, the largest degree of $q$ that
appears in $\sigma_{(\alpha\mid\beta)}\cdot\sigma_{\lambda}$ is~1. By
(\ref{thm_2_two}) we have         
\[
\sigma_{\lambda}*\sigma_{\mu}=\det(\sigma_{(\alpha_{i}\mid\beta_{j})})_{1\leq i,j\leq 
t}*\sigma_{\mu}=\sum_{\tau\in {\mathcal
S}_{n}}{\rm sign}(\tau)(\prod_{i=1}^{t}\sigma_{(\alpha_{i}\mid
\beta_{\tau(i)})})*\sigma_{\mu}.
\]
Thus $M\leq {\rm diag}(\lambda)$. Switching the roles of $\lambda$ and
$\mu$ gives $M\leq {\rm diag}(\mu)$.
\qed

It is an interesting problem to determine precisely which degrees of $q$ appear in
the product $\sigma_{\lambda}*\sigma_{\mu}$. 
Based on extensive computational evidence, we conjecture:

\begin{conjecture}
The product $\sigma_{\lambda}*\sigma_{\mu}$ involves $q^{d}$ for any 
$d\!\in\! [{\rm d_{min}},
{\rm d_{max}}]$.
\end{conjecture} 

        More seems to be true.

\begin{conjecture}
\label{second_conj}
Let $d\geq 1$ be an integer. If $\langle \lambda,\mu,\nu\rangle_{d}\neq 0$, then either
there exists $\alpha\supseteq \nu$ such that
$\langle \lambda,\mu,\alpha\rangle_{d-1}\neq 0$, or else 
$\langle \lambda,\mu,\alpha\rangle_{j}=0$ for all
$\alpha\subseteq l\times k$ and all $0 \leq j \leq d-1$.

\end{conjecture}
\section*{Acknowledgments}
This work was partially carried out while the author was visiting the
Fields Institute in Toronto, and later at the Isaac Newton Institute
in Cambridge. We thank both institutes for their hospitality. 
We are most deeply indebted to Sergey Fomin, 
whose guidance greatly improved the content and form of this paper. 
We benefited from conversations with Anders Buch, Bill Fulton, Arun Ram and Chris Woodward. 
Most of our computational investigations were done using Anders Buch's package~\cite{Buch}.

\end{document}